\documentclass[11pt]{article}
\usepackage{amsmath}
\usepackage{amssymb}
\usepackage{amsfonts}
\usepackage{mathrsfs}
\usepackage{shortvrb,psfrag}
\usepackage{enumerate}
\usepackage{color}
\usepackage[dvips]{graphicx}

\let\pa=\partial

\let\f=\frac
\let\p=\psi

\let\Om=\Omega

\def\R{\Bbb R}

\def\no{\noindent}
\def\na{\nabla}
\def\p{\partial}
\def\eqdef{\buildrel\hbox{\footnotesize def}\over =}
\def\endproof{\hphantom{MM}\hfill\llap{$\square$}\goodbreak}

\newcommand{\beq}{\begin{equation}}
\newcommand{\eeq}{\end{equation}}
\newcommand{\ben}{\begin{eqnarray}}
\newcommand{\een}{\end{eqnarray}}
\newcommand{\beno}{\begin{eqnarray*}}
\newcommand{\eeno}{\end{eqnarray*}}

\newtheorem{Theorem}{Theorem}[section]
\newtheorem{Definition}[Theorem]{Definition}
\newtheorem{Proposition}[Theorem]{Proposition}
\newtheorem{Lemma}[Theorem]{Lemma}

\newtheorem{Remark}[Theorem]{Remark}

\begin{document}
\title{Blow-up of critical norms for the 3-D Navier-Stokes equations}

\author{Wendong Wang$^\dag$ \,and\, Zhifei Zhang$^\ddag$\\[2mm]
{\small $^\dag$School of  Mathematical Sciences, Dalian University of Technology, Dalian 116024, P.R. China}\\[1mm]
{\small E-mail: wendong@dlut.edu.cn}\\[1mm]
{\small $^\ddag$School of  Mathematical Sciences, Peking University, Beijing 100871, P.R. China}\\[1mm]
{\small E-mail: zfzhang@math.pku.edu.cn}
}

\date{\today}
\maketitle

\begin{abstract}
Let $u=(u_h,u_3)$ be a smooth solution of the 3-D Navier-Stokes equations in $\R^3\times [0,T)$.
It was proved that if $u_3\in L^{\infty}(0,T;\dot{B}^{-1+3/p}_{p,q}(\R^3))$ for $3<p,q<\infty$ and $u_h\in L^{\infty}(0,T; BMO^{-1}(\R^3))$ with $u_h(T)\in VMO^{-1}(\R^3)$,
then $u$ can be extended beyond $T$. This result generalizes the recent result proved by
Gallagher, Koch and Planchon, which requires $u\in L^{\infty}(0,T;\dot{B}^{-1+3/p}_{p,q}(\R^3))$.
Our proof is based on a new interior regularity criterion in terms of one velocity component.
\end{abstract}

\setcounter{equation}{0}
\section{Introduction}
We consider the 3-D incompressible Navier-Stokes equations
\begin{equation}\label{eq:NS}
(NS)\left\{\begin{array}{l}
\partial_t u-\Delta u+u\cdot \nabla u+\nabla \pi=0,\\
{\rm div } u=0,\\
u|_{t=0}=u_0(x),
\end{array}\right.
\end{equation}
where $u(x,t)$ denotes the velocity of the fluid, and the scalar function $\pi(x,t)$ denotes the
 pressure.

The Navier-Stokes equations are scaling invariance in the sense:
if $(u,\pi)$ is a solution of (NS), then so does $(u_\lambda, \pi_\lambda)$, where
\ben\label{scaling}
u_\lambda=\lambda u(\lambda x, \lambda^2 t), \quad \pi_\lambda=\lambda^2\pi(\lambda x, \lambda^2 t).
\een
A Banach space $X$ is called a scaling invariant space with respect to the initial data if
$\|u_0\|_X=\|u_{0,\lambda}\|_X$. Some classical examples are $\dot H^\f12(\R^3), L^3(\R^3), \dot B^{-1+\f 3 p}_{p,q}(\R^3)$ and $BMO^{-1}(\R^3)$.
The local well-posedness of (NS) in $\dot H^\f12(\R^3), L^3(\R^3), \dot B^{-1+\f3 p}_{p,\infty}(\R^3)(3<p<\infty)$ and $BMO^{-1}(\R^3)$ was proved by Fujita-Kato \cite{FK}, Kato \cite{Kato}, Cannone-Meyer-Planchon \cite{CMP}  and Koch-Tataru \cite{KT} respectively. However, whether local smooth solution can be extended to a global one is an outstanding open problem in the mathematical fluid mechanics.

Ladyzhenskaya-Prodi-Serrin type criterion  states that if smooth solution $u$ satisfies
\beno
u\in L^q(0,T; L^p(\R^3))\quad \textrm{with} \quad \f 2 q+\f 3p\le 1\,\text{ for } p\ge 3,
\eeno
then $u$ can be extended after $t=T$, see \cite{Giga, Serrin, Struwe}  for example. The difficult endpoint case(i.e. $u\in L^\infty(0,T; L^3(\R^3))$ was solved by  Escauriaza-Seregin-\v{S}ver\'{a}k \cite{ESS}. This result in particular implies nonexistence of self-simliar type singularity \cite{NRS}.
\medskip

In this paper, we are concerned with the following interesting question:\medskip

{\bf If $u\in L^\infty(0,T;X)$ with $X$ a scaling invariant space, then $u$ can be extended beyond $t=T$?}

\medskip

The case of $X=L^3(\R^3)$ was proved by Escauriaza-Seregin-\v{S}ver\'{a}k(see \cite{KK, GKP1} for the alternative proofs by using the profile decomposition. The case of $X=\dot{B}^{-1+3/p}_{p,q}(\R^3)$ for $3<p<\infty, q<2p'$ was proved by Chemin-Planchon \cite{CP} by using self-improving bounds, and recently by Gallagher-Koch-Planchon \cite{GKP2} for the general case $3<p,q<\infty$. Recall that the following inclusion relations(see Lemma \ref{lem:inclusion})
\beno
\dot{H}^{1/2}(\R^3)\hookrightarrow L^3(\R^3)\hookrightarrow \dot{B}^{-1+3/p}_{p,q}(\R^3)(3<p<\infty)\hookrightarrow BMO^{-1}(\R^3).
\eeno
Thus, it is natural to ask whether the question holds for the case of $X= BMO^{-1}(\R^3)$.\medskip

In this paper, we prove a slightly weaker version, which requires that
the horizontal velocity $u_h\in L^{\infty}(0,T;BMO^{-1}(\R^3))$ with $u_h(T)\in VMO^{-1}(\R^3)$, while the vertical vertical $u_3\in L^{\infty}(0,T;\dot{B}^{-1+\frac3p}_{p,q}(\R^3))$.

\begin{Theorem}\label{thm:main}
Let $(u,\pi)$ be a smooth solution of (\ref{eq:NS}) in $\R^3\times [0,T)$. If $u$ satisfies
\beno
\|u_3\|_{L^{\infty}(0,T; \dot{B}^{-1+3/p}_{p,q})}+\|u_h\|_{L^{\infty}(0,T; BMO^{-1})}=M<\infty,
\eeno
for some $3<p,q<\infty$ and $u_h(x,T)\in VMO^{-1}(\R^3)$,
then $u$ can be extended after $t=T$.
\end{Theorem}

\begin{Remark}
Under the assumption of the theorem, it holds that
\beno
u(t)\in C_w\big([0,T];BMO^{-1}(\R^3)\big).
\eeno
Hence, $u(T)$ is well-defined as a function in $VMO^{-1}(\R^3)$. Thanks to the inclusion $\dot{B}^{-1+3/p}_{p,q}(\R^3)\hookrightarrow VMO^{-1}(\R^3)$,
our result improves the recent one proved by Gallagher, Koch and Planchon \cite{GKP2}, which requires $u\in L^{\infty}(0,T;\dot{B}^{-1+3/p}_{p,q}(\R^3))$.
\end{Remark}

The proof of \cite{GKP2} used the profile decomposition approach.
Our proof still follows blow-up analysis and backward uniqueness method developed by Escauriaza-Seregin-\v{S}ver\'{a}k \cite{ESS}.
A key ingredient is to prove a new interior regularity criterion in terms of one velocity component, which is independent of interest.
To state it, we introduce $G(f,p,q;r)\triangleq r^{1-\frac3p-\frac2q}\|f\|_{L^q_tL^{p}_x(Q_r)}$ and the scaling invariant quantities
\beno
C(u,r)\triangleq r^{-2}\int_{Q_r}|u|^3dxdt,\quad D(\pi,r)\triangleq r^{-2}\int_{Q_r}|\pi|^{\frac32}dxdt,
\eeno
where $Q_r=(-r^2,0)\times B_r(0)$.

\begin{Theorem}\label{thm:interior}
Let $(u,\pi)$ be a suitable weak solution of (\ref{eq:NS}) in $Q_1$. If $u$ satisfies
\beno
\sup_{0<r<1}G(u,p,q;r)\leq M<+\infty,
\eeno
for some $(p,q)$ with $1\le \frac3p+\frac2q<2$ and $1<q\leq \infty$,
then there exists a positive constant $\varepsilon$ depending on $p, q, M$
such that $(0,0)$ is a regular point if
\beno
G(u_3,p,q;r^*)\leq \varepsilon
\eeno
for some $r^*$ with $0<r^*<\min\{\frac1 2, (C(u,1)+D(\pi,1))^{-2}\}$.
\end{Theorem}

The second version is stated as follows.

\begin{Theorem}\label{thm:interior2}
Let $(u,\pi)$ be a suitable weak solution of (\ref{eq:NS}) in $Q_1$. If $(u,\pi)$ satisfies
\beno
C(u,1)+D(\pi,1)\leq M<+\infty,
\eeno
then there exists a positive constant $\varepsilon$ depending on $M$ such that if
\beno
\|u_3\|_{L^1(Q_1)}\leq \varepsilon,
\eeno
then $u$ is regular in $Q_{\f12}$.
\end{Theorem}

Compared with classical interior regularity criterions(see Proposition \ref{prop:small regularity-CKN} and Proposition \ref{prop:small regularity-GKT}),
new regularity criterion allows two components of the velocity to be large in
the local scaling invariant norm. The proof of Theorem \ref{thm:interior} and Theorem \ref{thm:interior2} used blow-up analysis. The key point is that under the assumptions of the above theorems the blow-up limit $v$ satisfies $v_3=0$ and
\beno
\left\{
\begin{array}{l}
\p_t v_h-\Delta v_h+v_h\cdot\na_h v_h+\na_h\pi=0,\\
\na_h\cdot v_h=0, \quad \p_{x_3}\pi=0,
\end{array}\right.
\eeno
which is similar to the situation considered by Cao-Titi \cite{CT1}.
Using the vorticity equation of the horizontal part $v_h$, it can be proved that the blow-up limit $v_h$ is regular(see Appendix for more details).

\section{Classical interior regularity criterion}

Let us first recall the definition of suitable weak solution introduced in \cite{CKN}.

\begin{Definition}
Let $\Om\subset \R^3$ and $T>0$. The pair $(u,\pi)$ is called a suitable weak solution
of (\ref{eq:NS}) in $Q_T=\Om\times (-T,0)$ if
\begin{enumerate}

\item $u\in L^{\infty}(-T,0;L^2(\Om))\cap L^2(-T,0;H^1(\Om))$ and $\pi\in L^{3/2}(Q_T)$;

\item $u$ satisfies (\ref{eq:NS}) in the sense of distribution;

\item $u$ satisfies the local energy inequality
\beno
&&\int_{\Om}|u(x,t)|^2\phi dx+2\int_{-T}^t\int_{\Om}|\nabla u|^2\phi dxds\\
&&\quad\leq
\int_{-T}^t\int_{\Om}|u|^2(\partial_s\phi+\triangle\phi)+u\cdot\nabla\phi(|u|^2+2\pi)dxds
\eeno
for a.e. $t\in [-T,0]$ and any nonnegative $\phi\in C_c^\infty(\R^3\times\R)$
vanishing in a neighborhood of the parabolic boundary of $Q_T$.

\end{enumerate}
\end{Definition}

Let $z_0=(x_0,t_0)\in \R^3\times \R$. We define a weak solution $u$ to be regular at $z_0$ if $u\in L^\infty(Q_r(z_0))$ for some $r>0$,
where $Q_r(z_0)=(-r^2+t_0,t_0)\times B_r(x_0)$.

The following quantities are invariant under the scaling (\ref{scaling}):
\beno
&&A(u,r,z_0)\triangleq\sup_{-r^2+t_0\leq t<t_0}r^{-1}\int_{B_r(x_0)}|u(y,t)|^2dy,\\
&&C(u,r,z_0)\triangleq r^{-2}\int_{Q_r(z_0)}|u(y,s)|^3dyds,\\
&&E(u,r,z_0)\triangleq r^{-1}\int_{Q_r(z_0)}|\nabla u(y,s)|^2dyds,\\
&&D(\pi,r,z_0)\triangleq r^{-2}\int_{Q_r(z_0)}|\pi(y,s)|^{\f32}dyds.
\eeno
We denote
\beno
&&G(f,p,q;r,z_0)\triangleq r^{1-\frac3p-\frac2q}\|f\|_{L^q_tL^{p}_x(Q_r(z_0))},\\
&&H(f,p,q;r,z_0)\triangleq r^{2-\frac3p-\frac2q}\|f\|_{L^q_tL^{p}_x(Q_r(z_0))},
\eeno
where the space-time norm $\|\cdot\|_{L^q_tL^{p}_x(Q_r(z_0))}$ is defined by
\beno
\|f\|_{L^q_tL^{p}_x(Q_r(z_0))}\eqdef\Big(\int_{t_0-r^2}^{t_0}\Big(\int_{B_r(x_0)}|f(x,t)|^pdx\Big)^{\f
q p}dt\Big)^\f 1q.
\eeno
For simplicity, we denote
$$A(u,r,(0,0))=A(u,r),\quad G(f,p,q;r,(0,0))=G(f,p,q;r)$$
and so on.\medskip

Now we recall the following $\varepsilon$-regularity results from \cite{CKN, GKT}.

\begin{Proposition}\label{prop:small regularity-CKN}
Let $(u,\pi)$ be a suitable weak solution of (\ref{eq:NS}) in $Q_1(z_0)$. There exists an
$\varepsilon_0>0$ such that if
\beno
\int_{Q_1(z_0)}|u(x,t)|^3+|\pi(x,t)|^{3/2}dxdt\leq \varepsilon_0,
\eeno
then $u$ is regular in $Q_{\f 12}(z_0)$. Moreover, $\pi$ can be replaced by
$\pi-(\pi)_{B_r}$ in the integral.
\end{Proposition}

\begin{Proposition}\label{prop:small regularity-GKT}
Let $(u,\pi)$ be a suitable weak solution of (\ref{eq:NS}) in $Q_1(z_0)$.
There exists an $\varepsilon_1>0$ such that if one of the following two conditions holds
\begin{enumerate}

\item $G(u,p,q;r,z_0)\leq \varepsilon_1$ for any $0<r<\f12$, where $1\leq \frac3p+\frac2q\leq 2$;\vspace{0.1cm}

\item $H(\na\times u,p,q;r,z_0)\leq \varepsilon_1$ for any $0<r<\f12$, where $2\leq \frac3{p}+\frac2{q}\leq 3$ and  $(p,q)\neq(1,\infty)$;
\end{enumerate}
then $u$ is  regular at the point $z_0$.
\end{Proposition}

Let us conclude this section by recalling the following lemmas from \cite{WZ}.

\begin{Lemma}\label{lem:pressure}
Let $(u,\pi)$ be a suitable weak solution of (\ref{eq:NS}) in $Q_1$ and $p>2, q\ge 2$. Then it holds that
\beno
H(\pi,p/2,q/2;r)\leq C\big(\frac{\rho}{r}\big)^{\frac{4}{q}+\frac6p-2}G(u,p,q;\rho)^2
+C\big(\frac{r}{\rho}\big)^{2-\frac{4}{q}}H(\pi,1,q/2;\rho)
\eeno
for any $0<4r<\rho<1$, where $C$ is a constant independent of $r,\rho$.
\end{Lemma}

\begin{Lemma}\label{lem:invariant}
Let $(u,\pi)$ be a suitable weak solution of (\ref{eq:NS}) in $Q_1$. If
\beno
&&\sup_{0<r<1}G(u,p,q;r)\leq M\quad\textrm{ with }\quad1\le\frac3p+\frac2q<2,\,1<q\leq \infty,
\eeno
then there holds
\beno
A(u,r)+E(u,r)+D(\pi,r)\leq C\big(r^{1/2}\big(C(u,1)+D(\pi,1)\big)+1\big)
\eeno
for any $0<r<1/2$, where $C$ is a constant depending on $p,q,M$.
\end{Lemma}

\section{BMO space and inclusion relations}

Recall that a local integrable function $f\in BMO(\R^3)$ if it satisfies
\beno
\sup_{x_0\in\R^3, R>0}\frac{1}{|B_R(x_0)|}\int_{B_R(x_0)}|f(x)-f_{B_R(x_0)}|dx<\infty.
\eeno
We denote by $VMO(\R^3)$ the complete space of $C_0^\infty(\R^3)$ under the norm of $BMO(\R^3)$. Thus, it holds that
\beno
\limsup_{R\downarrow 0}\frac{1}{|B_R(x_0)|}\int_{B_R(x_0)}|f(x)-f_{B_R(x_0)}|dx=0,\\
\limsup_{|x_0|\rightarrow +\infty}\frac{1}{|B_R(x_0)|}\int_{B_R(x_0)}|f(x)-f_{B_R(x_0)}|dx=0.
\eeno
We say that a function $u\in BMO^{-1}(\R^3)$ if there exist $U^j\in
BMO(\R^3)$ such that $u=\sum_{j=1}^3\p_jU^j$. $VMO^{-1}(\R^3)$ is
defined similarly. A remarkable property of $BMO$ function is
\ben\label{BMO-p1}
\sup_{x_0\in\R^3, R>0}\frac{1}{|B_R(x_0)|}\int_{B_R(x_0)}|f(x)-f_{B_R(x_0)}|^qdx<\infty.
\een
for any $1\le q<\infty$.

We also need the following Carleson measure characterization of BMO space.
We say that the tempered distribution $v$ is in $BMO(\R^3)$ if
\ben
\|v\|_{BMO}=\sup_{x_0\in \R^3,R>0}\Big( |B(x_0,R)|^{-1}\int_{B(x_0,R)}\int_0^{R^2}|\nabla e^{t\triangle}v|^2dtdx\Big)^{1/2}<\infty,
\een
and  $v$ is in $BMO^{-1}(\R^3)$ if
\ben\label{def:bmo-1}
\|v\|_{BMO^{-1}}=\sup_{x_0\in \R^3,R>0}\Big( |B(x_0,R)|^{-1}\int_{B(x_0,R)}\int_0^{R^2}|e^{t\triangle}v|^2dtdx \Big)^{1/2}<\infty.
\een
We refer to \cite{Graf} and \cite{Stein} for more introductions.
\medskip

For the homogeneous Besov space $\dot B^s_{p,q}(\R^3)$ for $s<0, p,q\in [1,\infty]$, we have the following characterization by
the heat flow:
\ben\label{def:besov}
\|f\|_{\dot{B}^s_{p,q}}\eqdef \big\|\|\tau^{-s/2} e^{\tau\triangle}f\|_{L^p(\R^3)}\big\|_{L^q(\R_+,\frac{d\tau}{\tau})}.
\een

Let us prove the following inclusion lemma.

\begin{Lemma}\label{lem:inclusion}
Let $p\in (3,\infty)$ and $q\in [1,\infty)$. It holds that
\beno
 \dot{B}^{-1+\frac3p}_{p,q}(\R^3)  \hookrightarrow VMO^{-1}(\R^3).
\eeno
\end{Lemma}

\no{\bf Proof.} Due to $\dot{B}^{-1+\frac3p}_{p,q_1}(\R^3)\subset \dot{B}^{-1+\frac3p}_{p,q_2}(\R^3))$ for $q_1\le q_2$, we may assume $q>2$. Let
\beno
I(x_0,R)=\Big(|B(x_0,R)|^{-1}\int_{B(x_0,R)}\int_0^{R^2}|e^{t\triangle}f|^2dtdx \Big)^{1/2},
\eeno
and $s=-1+\f3 p$. By H\"{o}lder inequality, we get
\beno
I(x_0,R)&\leq& C\Big[ R^{-\frac6p}\int_0^{R^2}\big(\int_{B(x_0,R)}|e^{t\triangle}f|^pdx\big)^{\frac2p} dt\Big]^{1/2}\\
&\leq& C\Big[ R^{-\frac6p}\int_0^{R^2}\big(\int_{B(x_0,R)}|e^{t\triangle}f|^pdx\big)^{\frac2p}t^{-\frac{q}{2}s-1} dt\Big]^{1/q}\Big[\int_0^{R^2}t^{(s+\frac2q)\frac{q}{q-2}}dt\Big]^{\frac{q-2}{2q}}\\
&\leq& C\Big[\int_0^{R^2}\big(\int_{B(x_0,R)}|e^{t\triangle}f|^pdx\big)^{\frac qp}t^{-\frac{q}{2}s-1} dt\Big]^{1/q}\\
&\leq& C\|f\|_{ \dot{B}^{-1+\frac3p}_{p,q}}
\eeno
which gives rise to $\dot{B}^{-1+\frac3p}_{p,q}(\R^3)\hookrightarrow BMO^{-1}(\R^3)$ by (\ref{def:bmo-1}) and (\ref{def:besov}).
Since $C_0^\infty(\R^3)$ is dense in $\dot{B}^{-1+\frac3p}_{p,q}(\R^3)$ for $p,q<\infty$, $f\in VMO^{-1}(\R^3)$.
\endproof

\begin{Lemma}\label{lem:interpolation}
If $u\in \dot B^{-1+3/p}_{p,\infty}(\R^3)\cap H^1(B_2(x_0))$ for $p\in (3,\infty)$ and some $x_0\in \R^3$, then we have
\beno
\|u\|_{L^2(B_1(x_0))}\le C(p)\|u\|_{\dot B^{-1+3/p}_{p,\infty}}^\theta\|u\|_{H^1(B_2(x_0))}^{1-\theta}\quad \text{for}\quad \theta=\f p {2p-3}.
\eeno
\end{Lemma}

\no{\bf Proof.}\,For any $\varphi(x)\in C_0^\infty(\R^3)$, we have
\beno
\|\varphi u\|_{\dot B^{-1+3/p}_{p,\infty}}\le C(\varphi)\|u\|_{\dot B^{-1+3/p}_{p,\infty}}.
\eeno
Thus, we may assume that $u\in H^1(\R^3)$ and it suffices to show that
\beno
\|u\|_{L^2(B_1(x_0))}\le C(p)\|u\|_{\dot B^{-1+3/p}_{p,\infty}}^\theta\|u\|_{\dot H^1}^{1-\theta}.
\eeno
We decompose $u$ into the low frequency part and high frequency part, i.e.,
\beno
u=u_N+u^N,\quad u_N={\cal F}^{-1}(\chi(\xi/2^N)\widehat u(\xi)),
\eeno
where $\chi(\xi)$ is a smooth cut-off function with $\chi(\xi)=1$ for $|\xi|\le 1$. For the low frequency, it can be proved by using the Littlewood-Paley
characterization of Besov space \cite{BCD} that
\beno
\|u_N\|_{L^2(B_1(x_0))}\le C\|u_N\|_{L^p(B_1(x_0))}\le C2^{N(1-3/p)}\|u\|_{\dot B^{-1+3/p}_{p,\infty}},
\eeno
while for the high frequency,
\beno
\|u^N\|_{L^2}\le C2^{-N}\|u\|_{\dot H^1}.
\eeno
Optimizing $N$ gives the lemma.\endproof

\section{Blow-up of the critical norm}

In this section, we prove Theorem \ref{thm:main} under the assumption of Theorem \ref{thm:interior} and Theorem \ref{thm:interior2}. By translation and scaling,
we may assume that $(u,\pi)$ is a smooth solution of (\ref{eq:NS}) in $\R^3\times [-1,0)$ satisfying the bound
\ben\label{ass:bound}
\|u_3\|_{L^{\infty}(-1,0;\dot{B}^{-1+3/p}_{p,q})}+\|u_h\|_{L^{\infty}(-1,0; BMO^{-1})}=M<\infty,
\een
and $u_h(x,0)\in VMO^{-1}(\R^3)$.
\medskip

Let us  first give a bound of local scaling invariant energy.

\begin{Lemma}\label{lem:energy}
Let $(u,\pi)$ be a smooth solution of (\ref{eq:NS}) in $\R^3\times [-1,0)$ satisfying (\ref{ass:bound}). Then
for $z_0\in\R^3\times(-\f14,0)$ and $0<r<\f12$, there holds
 \beno
A(u,r,z_0)+E(u,r,z_0)+D(\pi,r,z_0)\leq C\big(M, C(u,\f12,z_0), D(\pi,\f12,z_0)\big).
\eeno
\end{Lemma}

\no{\bf Proof.} Without loss of generality, we may assume $z_0=(0,0)$.
Let $\zeta(x,t)$ be a smooth function with $\zeta\equiv 1$ in $Q_r$
and $\zeta=0$ in $Q_{2r}^c$. Since $u\in
L^{\infty}(-1,0;BMO^{-1}(\R^3))$, there exists $U_i^j(x,t)\in
L^{\infty}(-1,0;BMO(\R^3))$ such that $u_i=\partial_j\cdot U^j_i$ for $i,j=1,2,3$.

We infer from H\"{o}lder inequality that
\beno
r^{-2}\int_{Q_{2r}}|u_i|^3\zeta^2dxdt&=&r^{-2}\int_{Q_{2r}}\sum_{j=1}^3\partial_j U_i^j\cdot u_i |u|\zeta^2dxdt\\
&\leq&Cr^{-2}\int_{Q_{2r}}|U-U_{B_{2r}}|(|\nabla u| |u|+|u|^2|\nabla\zeta|)dxdt,
\eeno
which gives rise to
\beno
C(u,r)&\leq& Cr^{-2}\Big(\int_{Q_{2r}}|U-U_{B_{2r}}|^6dxdt\Big)^{1/6}\Big(\int_{Q_{2r}}|\nabla u|^2dxdt\Big)^{1/2}\Big(\int_{Q_{2r}}|u|^3dxdt\Big)^{1/3}\\
&&+Cr^{-3}\Big(\int_{Q_{2r}}|U-U_{B_{2r}}|^3dxdt\Big)^{1/3}\Big(\int_{Q_{2r}}|u|^3dxdt\Big)^{2/3},
\eeno
from which and (\ref{BMO-p1}), it follows that
\beno
C(u,r)&\leq& C(M)\big(E(u,2r)^{1/2}C(u,2r)^{1/3}+ C(u,2r)^{2/3}\big)\\
&\leq& C(M)\big(1+E(u,2r)+C(u,2r)\big)^{5/6}.
\eeno
Using the local energy inequality, we can deduce that
\ben\label{eq:local energy inequ}
A(u,r)+E(u,r)\leq C\big(C(u,2r)^{2/3}+C(u,2r)+C(u,2r)^{1/3}D(\pi,2r)^{2/3}\big).
\een
This concludes that
\ben\label{eq:c(u,r)}
C(u,r)\leq C(M)\big(1+C(u,4r)+D(\pi,4r)\big)^{5/6}.
\een
Lemma \ref{lem:pressure} yields that for $0<4r<\rho<1$,
\ben\label{eq:pressure-d}
D(\pi,r)\leq C\big((\frac{r}{\rho})D(\pi,\rho)+(\frac{\rho}{r})^{2}C(u,\rho)\big).
\een

Let $F(r)=C(u,r)+D(\pi,r)$.  We infer from (\ref{eq:c(u,r)})--(\ref{eq:pressure-d})  that  for $0<16r<4\rho<1$,
\beno
F(r)&\leq& C(M)\big(1+C(u,4r)+D(\pi,4r)\big)^{5/6}+C(\frac{r}{\rho})D(\pi,\rho)\\
&&+C(M)(\frac{\rho}{r})^{2}\big(1+C(u,4\rho)+D(\pi,4\rho)\big)^{5/6}\\
&\leq& C(M)(\frac{r}{\rho})F(4\rho)+C(M)(\frac{\rho}{r})^{17}.
\eeno
Choosing $\theta$ such that $C(M)\theta<\f12$, we deduce that for any $0<r<1$,
\beno
F(\theta r)\le \f12F(r)+C(M,\theta).
\eeno
Then a standard iterative scheme  ensures that for $0<r<\f12$,
\beno
F(r)=C(u,r)+D(\pi,r)\leq C\big(M,C(u,\f12),D(\pi,\f12)\big),
\eeno
which along with (\ref{eq:local energy inequ}) implies the desired result.
\endproof
\medskip

\no{\bf Proof of Theorem \ref{thm:main}.}\,
Following \cite{ESS}, the proof is based on blow-up analysis and backward uniqueness of parabolic equations.
It suffices to prove that $u\in L^\infty\big(\R^3\times (-1,0)\big)$ by Ladyzhenskaya-Prodi-Serrin type criteria.
For this, we assume that $(0,0)$ is a singular point of $u$ without loss of generality.
\medskip

{\bf Step 1.} Blow-up analysis\vspace{0.1cm}

It follows from Lemma \ref{lem:energy} that for $z_0\in \R^3\times (-\f14,0)$ and $0<r<\f12$,
\ben\label{eq:bound of u_z0}
A(u,r,z_0)+E(u,r,z_0)+D(\pi,r,z_0)\leq C\big(M, C(u,1), D(\pi,1)\big).
\een
By the interpolation inequality, we get
\beno
C(u,r,z_0)\leq C\big(M, C(u,1), D(\pi,1)\big).
\eeno
If the point $(0,0)$ is not regular, Theorem \ref{thm:interior}  ensures that there exists a sequence of $r_k\downarrow 0$ such that
\ben\label{eq:sequence lower bdd}
r_k^{-2}\int_{Q_{r_k}}|u_3|^3dxdt\geq \varepsilon^3.
\een

Let
\beno
u^k(x,t)=r_ku(r_kx,r_k^2t),\quad \pi^k(x,t)=r_k^2\pi(r_kx,r_k^2t).
\eeno
For any $a>0, T>0$ and $z_0\in\R^3\times (-T,0)$, it follows from (\ref{eq:bound of u_z0}) that if $k$ large enough, we have
\ben\label{eq:bound of ukz_0}
A(u^k,a,z_0)+E(u^k,a,z_0)+C(u^k,a,z_0)+D(\pi^k,a,z_0)\leq C\big(M,C(u,1),D(\pi,1)\big).
\een
Then Aubin-Lions Lemma ensures that there exists $(v,\pi')$ such that for any $a,T>0$ (up to subsequence)
\beno
&&u^k\rightarrow v \quad {\rm in} \quad L^3(B_{a}\times (-T,0)),\\
&&u^k\rightarrow v \quad {\rm in} \quad C([-T,0];L^{9/8}(B_{a})),\\
&&\pi^k\rightharpoonup \pi' \quad {\rm in} \quad L^\f 32(B_{a}\times (-T,0)),
\eeno
as $k\rightarrow+\infty$. The lower semi-continuity of the norm gives that
\ben\label{eq:limit in Besov}
&&\|v_3\|_{L^{\infty}(-a^2,0; \dot{B}^{-1+3/p}_{p,q})}\leq M,\\\label{eq:limit bdd}
&&A(v,a,z_0)+E(v,a,z_0)+C(v,a,z_0)+D(\pi',a,z_0)\leq C\big(M,C(u,1),D(\pi,1)\big).
\een
Thanks to (\ref{eq:limit in Besov}), there exists a sequence $v_3^n\in C_0^\infty(\R^3\times (-a^2,0))$ so that for any $\ell\in [1,\infty)$,
\beno
\|v_3^n-v_3\|_{L^{\ell}(-T,0;\dot{B}^{-1+3/p}_{p,q})}\rightarrow0 \quad {\rm as}\,\, n\rightarrow\infty,
\eeno
and by (\ref{eq:limit bdd}), we have $\|v\|_{L^{{10}/{3}}(Q_a(z_0))}+\|\na v\|_{L^2(Q_a(z_0))}<\infty$. Thus, we have
\ben\label{eq:v3-infty}
\int_{Q_1(z_0)}|v_3|^3dxdt\rightarrow 0\quad {\rm as}\,\, |z_0|\rightarrow\infty.
\een
Indeed, for any $\epsilon>0$, we first take $n$ large enough so that
\beno
\|v_3^n-v_3\|_{L^{\ell}(-T,0;\dot{B}^{-1+3/p}_{p,q})}<\epsilon.
\eeno
Then take $|z_0|$ big enough so that $v_3^n=0$ in $Q_2(z_0)$. By Lemma \ref{lem:interpolation}, we have
\beno
\int_{Q_1(z_0)}|v_3|^2dxdt&=&\int_{Q_1(z_0)}|v_3^n-v_3|^2dxdt\\
&\le& C\|v_3^n-v_3\|_{L^{\ell}(-T,0;\dot{B}^{-1+3/p}_{p,q})}^{2\theta}
\|v_3\|_{L^2_tH^1_x(Q_2(z_0))}^{2(1-\theta)}\le C\epsilon^{2\theta},
\eeno
which implies (\ref{eq:v3-infty}) by the interpolation.

Thanks to (\ref{eq:limit bdd}) and (\ref{eq:v3-infty}), Theorem \ref{thm:interior2} ensures that there exists $R>0$ such that $v$ is regular
in $(\R^3\backslash{B_R})\times (-T,0)$ with the bound
\ben\label{eq:bound of v}
|v(x,t)|+|\nabla v(x,t)|\leq C \quad\text{for}\,\, (t,x)\in (\R^3\backslash{B_R})\times (-T,0).
\een

{\bf Step 2.} Prove $v(x,0)=0$ in $\R^3$\medskip

Due to  $u(x,0)\in VMO^{-1}(\R^3)$, $u=\nabla\cdot U$ with $U(x,0)\in VMO(\R^3)$.
Then for any $\varphi\in C_c^\infty(B_a)$,
\beno
&&\Big|\int_{B_a}v(x,0)\varphi dx\Big|\\
&&\leq C\int_{B_a}|v(x,0)-u^k(x,0)|dx+\big|\int_{B_a}u^k(x,0)\varphi dx\big|\\
&&\leq C\int_{B_a}|v(x,0)-u^k(x,0)|dx+r_k^{-3}\big|\int_{B_{ar_k}}U(y,0)(\na\varphi)(y/R_k)dy\big|\\
&&\leq C\int_{B_a}|v(x,0)-u^k(x,0)|dx+Cr_k^{-3}\int_{B_{ar_k}}|U(y,0)-U_{B_{ar_k}}|dy\\
&&\longrightarrow 0\quad\textrm{ as } k\rightarrow\infty.
\eeno
This shows that $v(x,0)=0.$\medskip

{\bf Step 3.} Backward uniqueness and unique continuation argument\vspace{0.1cm}

Let $w=\nabla\times v$.  By (\ref{eq:bound of v}) and $v(x,0)=0$, we get $w(x,0)=0$ and
$$
|\partial_t w-\triangle w|\leq C(|w|+|\nabla w|)\quad \textrm{in}\quad (\R^3\backslash{B_R})\times (-T,0).
$$
By the backward uniqueness of parabolic operator \cite{ESS}, we deduce that
$w=0$ in $(\R^3\backslash{B_R})\times (-T,0).$ Similar arguments as in \cite{ESS}, by using spacial unique continuation,
give $w=0$ in $\R^3\times (-T,0)$. Hence, we get
\beno
\triangle v\equiv 0\quad  \textrm{in}\quad  \R^3\times (-T,0).
\eeno
This implies $v_3=0$ in $\R^3\times (-T,0)$, since $v_3(\cdot,t)\in \dot{B}^{-1+3/p}_{p,q}(\R^3)$.
This leads to an obvious contradiction with the fact (\ref{eq:sequence lower bdd}).
Thus, $(0,0)$ is a regular point.
\endproof

\section{New interior regularity criterion}

In this section, we first prove two new interior regularity criterion in terms of one velocity component. Then we present
an applications to type I singularity.

\subsection{Interior regularity criterion}

In this subsection, we prove Theorem \ref{thm:interior} and Theorem \ref{thm:interior2} under the assumption of Proposition \ref{prop:blowup}.\medskip

\no{\textbf{Proof of Theorem \ref{thm:interior}}. The proof uses blow-up argument. Assume that the statement of the theorem is false.
Then there exist constants $p, q, M$ and a sequence of suitable weak solutions $(u^k,\pi^k)$ of (\ref{eq:NS}) in $Q_1$ ,
which are singular at $(0,0)$ and satisfy
\beno
&&G(u^k,p,q;r)\leq M \quad {\rm for\,\, all}\quad 0<r<1,\\
&&G(u_3^k,p,q;r_k)\leq \frac1k,
\eeno
where
$$0<r_k<\min\big\{\frac12, (C(u^k,1)+D(\pi^k,1))^{-2}\big\}.$$
Then it follows from Lemma \ref{lem:invariant} that
\beno
A(u^k,r)+E(u^k,r)+D(\pi^k,r)\leq C_0(M,p,q)
\eeno
for any $0<r<r_k.$ Set
$$v^k(x,t)=r_ku^k(r_kx,r_k^2t), q^k(x,t)=r_k^2\pi^k(r_kx,r_k^2t).$$
Thus, it holds that
\beno
A(v^k,r)+E(v^k,r)+D(\pi^k,r)\leq C_0(M,p,q)
\eeno
for any $0<r<1$. Lions-Aubin's lemma ensures that there exists a suitable weak solution $({v}, {\pi})$ of (\ref{eq:NS})
such that (up to subsequence),
\beno
&&v^k\rightarrow {v}\quad {\rm in} \quad L^3(Q_\f12),\quad \pi^k\rightharpoonup \pi\quad {\rm in} \quad L^\f {3}2(Q_\f12),\\
&& v_3^k\rightharpoonup 0\quad {\rm in} \quad L^q((-\f14,0); L^p(B_\f12)),
\eeno
as $k\rightarrow+\infty$.
So, the blow-up limit $v$ satisfies $v_3=0, \partial_3 \pi=0$ and
\beno
\partial_tv_h+v_h\cdot\nabla_h v_h+\nabla_h \pi-\triangle v_h=0, \quad \nabla_h\cdot v_h=0.
\eeno
The limit $v_h$ is bounded in $Q_\f12$ by Proposition \ref{prop:blowup}, which will contradicts with the fact that $(0,0)$ is a singular point of $v^k$.
Indeed by Proposition \ref{prop:small regularity-CKN} and Lemma \ref{lem:pressure}, for any $0<r<1/4$,
\beno
\varepsilon_0 &\leq&\liminf_{k\rightarrow\infty}r^{-2}\int_{Q_r}|v^k|^3+|\pi^k|^{3/2}dxdt\\
&\leq& C\liminf_{k\rightarrow\infty}\Big(C({v},r)+\frac{r}{\rho}D(\pi^k,\rho)+(\frac{\rho}{r})^2C(v^k,\rho)\Big)\\
&\leq& C\Big(C({v},r)+\frac{r}{\rho}+(\frac{\rho}{r})^2C({v},\rho)\Big)\\
&\leq& Cr^{1/2}\quad (\textrm{by choosing}\quad \rho=r^{1/2}),
\eeno
which is a contradiction if we take $r$ small enough.\endproof

\medskip

To prove Theorem \ref{thm:interior2}, we need the following lemma from \cite{WZ}.

\begin{Lemma}\label{lem:c}
Let $(u,\pi)$ be a suitable weak solution of (\ref{eq:NS}) in $Q_1$ and ${D}(\pi,1)\leq M$.
Then $u$ is regular at $(0,0)$ if
\beno
C(u,r_0)\leq c\varepsilon_0^{9/5}r_0^{8/5}\quad\textrm{ for some }0<r_0\leq 1.
\eeno
Here $c$ is a small constant depending on $M$.
\end{Lemma}

\no{\textbf{Proof of Theorem \ref{thm:interior}}. We use the contradiction argument.
Assume that the statement of the theorem is false. Then there exist $M$ and a sequence of suitable weak solutions $(u^k,\pi^k)$ of (\ref{eq:NS}) in $Q_1$,
which are singular at a point $z_0\in Q_\f12$(we assume $z_0=(0,0)$) and satisfy
\beno
C(u^k,1)+D(\pi^k,1)\leq M,\quad \|u_3^k\|_{L^1(Q_1)}\leq \frac 1k.
\eeno
Then by the local energy inequality, it is easy to show that
\beno
A(u^k,3/4)+E(u^k,3/4)\leq C(M).
\eeno
By using Lions-Aubin's lemma, there exists a suitable weak solution $(v, \pi')$ of (\ref{eq:NS}) such that (up to subsequence),
\beno
u^k\rightarrow v,\quad  u_3^k\rightarrow 0  \quad {\rm in} \quad L^3(Q_\f12),\quad
\pi^k\rightharpoonup \pi'\quad {\rm in} \quad L^\f {3}2(Q_\f12),
\eeno
as $k\rightarrow+\infty$. Note that $v_3=0$ implies $\partial_3 \pi'=0$, hence
\beno
\partial_tv_h+v_h\cdot\nabla_h v_h+\nabla_h \pi'-\triangle v_h=0, \quad \nabla_h\cdot v_h=0.
\eeno
The limit $v_h$ is bounded in $Q_\f12$ by Proposition \ref{prop:blowup}.
Since $(0,0)$ is a singular point of $u^k$,  we have by Lemma \ref{lem:c} that for any $0<r<1/4$,
\beno
c\varepsilon_0^{9/5}r^{8/5}&\leq& \lim_{k\rightarrow\infty}r^{-2}\int_{Q_r}|u^k|^3dxdt\\
&\leq& \lim_{k\rightarrow\infty} C(v,r)
\leq C(M)r^{3},
\eeno
which is a contradiction by letting $r\rightarrow 0$. \endproof

\subsection{Application to type I singularity}

In two important works \cite{CSTY, KNSS}, the authors excluded the existence of type I singularity\big(i.e., $\|u(t)\|_{L^\infty}\le \f C {\sqrt{-t}}$\big) for the axisymmetric Navier-Stokes equations. The general case remains open. As an application of interior regularity criterion, we prove nonexistence
of type I singularity in the case when one velocity component has a better control.

\begin{Theorem}
Let $(u,\pi)$ be a suitable weak solution of (\ref{eq:NS}) in $Q_1$ and satisfy
\beno
|u_h(x,t)|\leq \frac{M}{\sqrt{-t}}.
\eeno
If $u_3\in L^q_tL^p_x(Q_1)$ with $\frac3p+\frac2q\le 1$ for $3<p\leq \infty$ or there exists $\epsilon>0$ depending on $M$ so that
$|u_3(x,t)|\leq \frac{\epsilon}{\sqrt{-t}}$, then $u$ is regular at $(0,0)$.
\end{Theorem}

\no{\bf Proof.}\,In the case when $u_3\in L^q_tL^p_x(Q_1)$, Lemma \ref{lem:invariant} ensures that there exists $r_0>0$  so that
for any $0<r<r_0$,
\beno
&&C(u,r)+D(\pi,r)\leq M,\quad \|u_3\|_{L^q_tL^p_x(Q_r)}\le \varepsilon,
\eeno
which implies that $u$ is regular at $(0,0)$ by Theorem \ref{thm:interior}.
The case of $|u_3(x,t)|\leq \frac{\epsilon}{\sqrt{-t}}$ follows from Theorem \ref{thm:interior2} with $p=\infty$ and $q=\f32$.\endproof
\medskip

Recently, there are many interesting works devoted to the regularity criterions involving one velocity component,
see \cite{CT2,CT3,CZ, CZZ, KZ,PZ} and references therein. However, the obtained regularity conditions are not scaling invariant
in these works except \cite{CZ,CZZ}, where they showed that if $u_3\in L^p(0,T;H^{\f12+\f2p}(\R^3))$ for $p\in (4,\infty)$, then the solution can be extended after $t=T$. The question of whether the solution is regular in the class $u_3\in L^q(0,T;L^p(\R^3))$ with $\frac3p+\frac2q=1$ for $3<p\leq \infty$ is still open.

\section{Appendix}

In this appendix, we study the regularity of the blow-up limit introduced in the proof of interior regularity criterion, which satisfies
\begin{equation}\label{eq:blow}
\left\{\begin{array}{l}
\p_t v_h-\Delta v_h+v_h\cdot\na_h v_h+\na_h\pi=0,\\
\na_h\cdot v_h=0,\quad \p_{x_3}\pi=0.
\end{array}\right.
\end{equation}
Here $\na_h=(\pa_{x_1}, \pa_{x_2})$.

\begin{Proposition}\label{prop:blowup}
Let  $(v_h, \pi)$ be a suitable weak solution of (\ref{eq:blow}) in $Q_1$.  Assume that
\beno
\|v_h\|_{L^\infty_tL^2_x(Q_1)}+\|\na v_h\|_{L^2(Q_1)}\le M.
\eeno
Then it is regular in $Q_\f14$ and $\|v_h\|_{L^\infty(Q_\f14)}\le C(M)$.
\end{Proposition}

\no{\bf Proof.}\,
By Proposition \ref{prop:small regularity-GKT}, we only need to  prove that $\|\nabla v_h\|_{L^{\infty}_tL^2_x(Q_{\frac13})}$ is finite.
 Using the fact that  the $\frac12$ dimensional Hausdorff measure of the set of singular time is zero (see \cite{CKN}), without loss of generality we may assume $v_h\in C^{\infty}(Q_1)$.  Let $w_h=\partial_1v_2-\partial_2v_1$ be the vorticity of $v_h$. Then $w_h$ satisfies
\begin{equation}\label{eq:blow-wh}
\p_t w_h-\Delta w_h+v_h\cdot\na_h w_h=0\quad {\rm in}\quad Q_1.
\end{equation}
While, $d=\partial_3v_h$ satisfies
\begin{equation}\label{eq:blow-3 vh}
\p_td-\Delta d+v_h\cdot\na_h d+d\cdot\na_hv_h=0\quad {\rm in}\quad Q_1.
\end{equation}

{\bf Step 1}. Estimate of $\|w_h\|_{L^{\infty}_tL^2_x(Q_{\frac34})}$\medskip

Let $\zeta$ be a cutoff function, which vanishes
outside of $Q_{1}$ and equals 1 in $Q_{\frac34}$, and satisfies
$$|\nabla\zeta|+ |\partial_t\zeta|+|\nabla^2\zeta|\leq C.$$
Acting on $2u\zeta^4$ on both sides of (\ref{eq:blow-wh}), we obtain
\beno
I&\triangleq&\sup_t\int_{B_{1}}|w_h|^2\zeta^4 dx+2\int_{Q_{1}}|\nabla (w_h\zeta^2)|^2 dxdt\\
&\leq& \int_{Q_{1}}\big[\partial_t(\zeta^4)+2|\nabla(\zeta^2)|^2\big]|w_h|^2dxdt+\int_{Q_{1}}v_h\cdot\nabla(\zeta^4)|w_h|^2dxdt\\
&\triangleq&I_1+I_2.
\eeno
Obviously, $I_1\leq C(M)$. For $I_2$, we get by Sobolev inequality that
\beno
I_2&\leq& C\int_{Q_{1}}|v_h||w_h|^{\frac12} | (w_h\zeta^2)|^{\frac32}dxdt\\
&\leq& C\|v_h\|_{L^{\infty}_tL^2_x(Q_1)}\|w_h\|_{L^{2}(Q_1)}^{\frac12}\|w_h\zeta^2\|_{L^{2}_tL^6_x(Q_1)}^{\frac32}\\
&\leq& C(M)I^{\frac34}.
\eeno
This shows that
\beno
I\leq  C(M)+C(M)I^{\frac34},
\eeno
which implies that $\|w_h\|_{L^{\infty}_tL^2_x(Q_{\frac34})}<\infty.$\medskip

{\bf Step 2}. Estimate of $\|\nabla_hv_h\|_{L^{\infty}_tL^2_x(Q_{\frac12})}$\medskip

Since $w_h$ is the vorticity of $v_h$ and $v_h$ is divergence-free, the elliptic estimate gives
\beno
\|\nabla_hv_h\|_{L^2_x(B_{\frac12})}\leq C\big( \|w_h\|_{L^2_x(B_{\frac34})}+  \|v_h\|_{L^2_x(B_{\frac34})}\big)\le C(M).
\eeno

{\bf Step 3.} Estimate of $\|d\|_{L^{\infty}_tL^2_x(Q_{\frac12})}$}.

Let $\eta$ be a cutoff function, which vanishes
outside of $Q_{\frac12}$ and equals 1 in $Q_{\frac13}$, and satisfies
$$|\nabla\eta|+ |\partial_t\eta|+|\nabla^2\eta|\leq C_0.$$
Acting $2d\eta^4$ on both sides of (\ref{eq:blow-3 vh}),
 we obtain
\beno
II&\triangleq&\sup_t\int_{B_{1}}|d|^2\eta^4 dx+2\int_{Q_{1}}|\nabla (d\eta^2)|^2 dxdt\\
&\leq& \int_{Q_{1}}\big[\partial_t(\eta^4)+2|\nabla(\eta^2)|^2\big]|d|^2dxdt
+\int_{Q_{1}}v_h\cdot\nabla(\eta^4)|d|^2dxdt\\
&&- \int_{Q_{1}}d\cdot\na_hv_h\cdot \big(2d\eta^4\big)dxdt\\
&\triangleq& II_1+II_2+II_3.
\eeno
Obviously, $II_1\leq C(M)$. Similar to $I_2$, we have
\beno
II_2\leq C(M)II^{\frac34}.
\eeno
Thanks to  $\|\nabla_hv_h\|_{L^{\infty}_tL^2_x(Q_{\frac12})}\le C(M)$, we have
\beno
II_3
&\leq& C\|\nabla_hv_h\|_{L^{\infty}_tL^2_x(Q_{\frac12})}\|d\|_{L^{2}(Q_1)}^{\frac12}\|d\eta^2\|_{L^{2}_tL^6_x(Q_1)}^{\frac32}\\
&\leq& C(M)II^{\frac34}.
\eeno
This shows that
\beno
II\leq  C(M)+C(M)II^{\frac34},
\eeno
which implies that $\|d\|_{L^{\infty}_tL^2_x(Q_{\frac12})}\le C(M).$\endproof

\vspace{0.2cm}

\bigskip

\noindent {\bf Acknowledgments.}
W. Wang is supported NSF of China 11301048 and the Fundamental Research Funds for the Central Universities.
Z. Zhang was partially supported by NSF of China under Grant 11371037 and 11425103.

\end{document}